\documentclass{article}

\usepackage[final, nonatbib]{nips_2016}

\usepackage[numbers]{natbib}

\usepackage[utf8]{inputenc} 
\usepackage[T1]{fontenc}    
\usepackage{hyperref}       
\usepackage{url}            
\usepackage{booktabs}       
\usepackage{amsfonts}       
\usepackage{nicefrac}       
\usepackage{microtype}      

\usepackage{graphicx}
\usepackage{amsmath}
\usepackage{amssymb}
\usepackage{xspace}
\usepackage[shortcuts,cyremdash]{extdash}
\graphicspath{{./figures/}}

\usepackage{nth}
\usepackage{algorithm}
\usepackage[noend]{algpseudocode}
\usepackage{amsthm}

\def\EE{{\mathbb E}}

\def\RR{{\mathbb R}}

\def\x{{\mathbf x}}

\def\A{{\mathbf A}}
\def\B{{\mathbf B}}
\def\C{{\mathbf C}}
\def\G{{\mathbf G}}

\def\D{{\mathbf D}}

\def\I{{\mathbf I}}

\def\M{{\mathbf M}}

\def\X{{\mathbf X}}

\def\d{{\mathbf d}}

\def\x{{\mathbf x}}

\newcommand{\balpha}{\boldsymbol{\alpha}}
\newcommand{\bbeta}{\boldsymbol{\beta}}

\renewcommand\P{{\mathbf P}}

\def\trace{{\mathrm{Tr}\;}}

\newcommand{\argmin}{\operatornamewithlimits{\mathrm{argmin}}}

\newcommand{\equaldef}{\stackrel{\textrm{def}}{=}}

\newcommand{\vertiii}[1]{{\left\vert\kern-0.25ex\left\vert\kern-0.25ex\left\vert #1
    \right\vert\kern-0.25ex\right\vert\kern-0.25ex\right\vert}}

\newtheoremstyle{ieee}
  {}
  {}
  {}
  {}
  {\itshape}
  {.}
  { }
  {\thmname{#1}\thmnumber{ #2}\thmnote{ (#3)}}%

\newtheoremstyle{assumption}
{}
{}
{}
{1em}
{\itshape}
{}
{ }
{\thmname{#1} \thmnumber{#2} \textnormal{\textit{\thmnote{#3.}}}}%

\theoremstyle{ieee}

\theoremstyle{assumption}

\theoremstyle{ieee}

\theoremstyle{ieee}

\theoremstyle{ieee}

\theoremstyle{ieee}

\algnewcommand\algorithmicswitch{\textbf{switch}}
\algnewcommand\algorithmiccase{\textbf{case}}
\algnewcommand\algorithmicassert{\texttt{assert}}
\algnewcommand\Assert[1]{\State \algorithmicassert(#1)}%
\algdef{SE}[SWITCH]{Switch}{EndSwitch}[1]{switch #1}{\algorithmicend\ \algorithmicswitch}%
\algdef{SE}[SWITCH]{Block}{EndBlock}[1]{#1}{}%
\algdef{SE}[CASE]{Case}{EndCase}[1]{\algorithmiccase\ #1}{\algorithmicend\ \algorithmiccase}%
\algtext*{EndSwitch}%
\algtext*{EndCase}%

\algnewcommand{\LineComment}[1]{\State \(\triangleright\) #1}
\algnewcommand{\Input}[1]{\State \textbf{Input}: #1}
\algnewcommand{\Output}[1]{\State \textbf{Output}: #1}

\newcommand{\somf}{\textsc{somf}\xspace}

\providecommand{\tightlist}{%
  \setlength{\itemsep}{0pt}\setlength{\parskip}{0pt}}

\title{Subsampled online matrix factorization\\ with convergence guarantees}

\author{
  Arthur Mensch\\
  Inria Parietal\\
  Saclay, France\\
  \And
  Julien Mairal\\
  Inria Thoth\\
  Grenoble, France\\
  \And
  Ga\"el Varoquaux\\
  Inria Parietal\\
  Saclay, France\\
  \And
  Bertrand Thirion\\
  Inria Parietal\\
  Saclay, France
}
\begin{document}
\maketitle
\vspace*{-3.4em}
\begin{center}
    \texttt{firstname.lastname@inria.fr}
\end{center}
\vspace*{2.4em}

\begin{abstract}

We present a matrix factorization algorithm that scales to input
matrices that are large in both dimensions (i.e., that contains more
than 1TB of data).  The algorithm streams the matrix columns while
subsampling them, resulting in low complexity per iteration and
reasonable memory footprint. In contrast to previous online matrix
factorization methods, our approach relies on low-dimensional
statistics from past iterates to control the extra variance introduced
by subsampling. We present a convergence analysis that guarantees us
to reach a stationary point of the problem. Large speed-ups can be
obtained compared to previous online algorithms that do not perform
subsampling, thanks to the feature redundancy that often exists in
high-dimensional settings.
\end{abstract}

\paragraph{Setup.}The goal of matrix factorization is to
decompose a matrix $\X \in \RR^{p \times n}$ -- typically $n$ signals of
dimension $p$ -- as a product of two smaller matrices:
\begin{equation}
    \X \approx \D \A
 \quad \text{with}\quad\D \in \RR^{p \times k}, \;\A \in \RR^{k \times n},
\end{equation}
with potential sparsity or structure requirements on $\D$ and~$\A$.
We consider a sample stream $(\x_t)_{t \geq 0}$ that cycles into the columns ${\{\x^{(i)}\}}_i$
of $\X$.
Matrix factorization can be formulated as a non-convex optimization
problem, where the factor $\D$ (the \textit{dictionary}) minimizes the following empirical risk:
\begin{equation}
    \label{eq:empirical-risk}
    \D = \argmin_{\D \in \mathcal{C}}\;
    \bar f \equaldef \frac{1}{n} \sum_{i=1}^n
    f^{(i)}(\D),\quad\text{where}\quad
    f^{(i)}(\D) =
      \min_{\balpha \in \RR^k} \frac{1}{2}
     \bigl\|
     \x^{(i)}
     - \D \balpha
     \bigr\|_2^2 + \lambda \, \Omega(\balpha).
\end{equation}
$\mathcal{C}$ is a column-wise separable convex set of $\RR^{p \times
k}$, and $\Omega : \RR^p \rightarrow \RR$ is a penalty over the code. The
problem of \textit{dictionary learning} \cite{olshausen_sparse_1997, agarwal_learning_2014} sets $\mathcal{C} = \mathcal{B}_2^k$
and $\Omega = \Vert \cdot \Vert_1$. Due to the sparsifying effect of $\ell_1$
penalty \cite{tibshirani_regression_1996}, the dictionary forms a basis in which the data admit a
\textit{sparse} representation. Setting $\mathcal{C} = \mathcal{B}_ 1^k$
and $\Omega = \Vert \cdot \Vert_2^2$ yields a data-adapted \textit{sparse basis}, akin to sparse PCA\cite{zou_sparse_2006}. The
algorithm presented here accommodates \textit{elastic-net} penalties
$\Omega(\balpha) \triangleq (1- \nu) \Vert \balpha \Vert_2^2 + \nu \Vert \balpha \Vert_1$,
and elastic-net ball-constraints $\mathcal{C} \triangleq \{ \D \in \RR^{p \times k}/\:
 \Vert \d^{(j)} \Vert \triangleq \Vert \d^{(j)} \Vert_1 + (1 - \mu) \Vert \d^{(j)} \Vert_2^2 \leq 1 \}$.
\paragraph{Problem.}For many applications of matrix factorization, datasets are
growing in both sample number and sample dimension. The online matrix
factorization algorithm \cite{mairal_online_2010} can handle large numbers of
samples but was designed to work in relatively small dimension. Recent work
\cite{mensch_dictionary_2016} has adapted this algorithm to handle very high
dimensional dataset. Though it demonstrates good empirical performance, the
proposed algorithm yields sequences of
iterates with non vanishing variance and is not asymptotically convergent.

\paragraph{Contribution.}We address this issue and correct some aspects of
the algorithm in \cite{mensch_dictionary_2016}  to establish
convergence and correctness. We thus introduce a new method to efficient solve
\eqref{eq:empirical-risk} for numerous high-dimensional data --- large $p$,
large $n$ --- with theoretical guarantees.

As in \cite{mensch_dictionary_2016}, we perform \textit{subsampling}
at each iteration of the online matrix factorization algorithm. The sample stream is
downsampled in a \textit{stochastic} manner ---we observe different subsets of successive colunms.
We thus each step of the online algorithm in a space
of reduced dimension $q < p$. Unlike \cite{mensch_dictionary_2016},
we control the variance introduced by the subsampling to obtain
convergence. For this, we rely on low-dimensional statistics kept from the past, as many recent stochastic algorithms
\cite{schmidt_minimizing_2013, johnson_accelerating_2013, defazio_saga:_2014}.\vspace*{-.4em}

\section{Algorithm}
\paragraph{Original online algorithm.}The problem \eqref{eq:empirical-risk} can be solved online, following \cite{mairal_online_2010}.
At each iteration $t$, a sample $\x_t$ is drawn from one of the columns $\{ \x^{(i)} \}_{1 \leq i \leq n}$ of $\X$.
 Its code $\balpha_t$ is computed from the previous dictionary $\D_{t-1}$: $\balpha_s \triangleq \argmin_{\balpha \in \RR^p} \frac{1}{2}
 \bigl\|
 \x_s
 - \D_{s-1} \balpha
 \bigr\|_2^2 + \lambda \, \Omega(\balpha)$. Then, $\D_t$ is updated as
 \begin{equation}
     \label{eq:minimization}
     \D_t \in \argmin_{\D \in \mathcal{C}} \bar g_t(\D)
     \triangleq \Big(
     \frac{1}{t}\sum_{s=1}^t \frac{1}{2}
     \bigl\|
     \x_s
     - \D \balpha_s
     \bigr\|_2^2 + \lambda \Omega(\balpha_s) \Big),
 \end{equation}
In other words, $\D_t$ is chosen to be the best dictionary that relates past
codes $(\balpha_s)_{s \leq t}$ to past samples $(\x_t)_{s \leq t}$. Codes are not recomputed from the current iterate $\D$,
which would be necessary to compute the true past loss fonction $\bar f_t(\D)$, of which $\bar g_t$ is a strongly-convex upper-bound:
\begin{align}
    \label{eq:bar_ft}
    \bar f_t(\D) &\triangleq \frac{1}{t} \sum_{s=1}^t \min_{\balpha \in \RR^p} \frac{1}{2}
        \bigl\|
        \x_s
        - \D \balpha
        \bigr\|_2^2\!+ \lambda \Omega(\balpha) \leq \bar g_t(\D)
\end{align}
It can
be shown \cite[see][for theoretical grounding]{mairal_stochastic_2013} that
minimizing $(\bar g_t)_t$ yields a sequence of iterates that is asymptotically a critical point
of $\bar f$ defined in~\eqref{eq:empirical-risk}. $\bar g_t$ can be minimized efficiently by projected block
coordinate descent, which makes it useful in practice. Indeed, minimizing $\bar g_t$
is equivalent to minimizing the quadratic function
\begin{equation}
    \label{eq:full_quadratic}
    \D \to\frac{1}{2} \trace (\D^\top \D \bar \C_t)
      - \trace (\D^\top \B_t),\:\text{where}\quad
      \bar \B_t = \frac{1}{t} \sum_{s=1}^t \x_s \balpha_s^\top, \quad
      \bar \C_t = \frac{1}{t} \sum_{s=1}^t \balpha_s \balpha_s^\top.
\end{equation}
Its gradient $\nabla \bar g_t: \D \to \D \bar \C_t - \bar \B_t$ can be tracked
online by updating $\bar \C_t$ and $\bar \B_t$ at each iteration:
\begin{equation}
    \label{eq:parameter-aggregation}
    \bar \C_t = (1 - \frac{1}{t}) \bar \C_{t-1}
    + \frac{1}{t} \balpha_t \balpha_t^\top \qquad
    \bar \B_t = (1 - \frac{1}{t}) \bar \B_{t-1}
    + \frac{1}{t} \x_t \balpha_t^\top.
\end{equation}
Those two statistics are thus sufficient to yield the sequence $(\D_t)_t$.
The weight $\frac{1}{t}$ used above can be replaced by a more general
$w_t$. In addition, the online algorithm has a minibatch extension.

\paragraph{Subsampled algorithm.}We adapt the algorithm from~\cite{mairal_online_2010} to
handle large sample dimension $p$. The complexity of this algorithm
linearly depends on the dimension $p$ in three aspects:
\tightlist
\begin{itemize}
    \item $\x_t \in \RR^p$ is
    used to compute the code $\balpha_t$,
    \item it is used to update the surrogate parameters
    $\bar \C_t \in \RR^{p\times k}$,
    \item $\D_t \in \RR^{p\times k}$ is fully updated at each iteration.
\end{itemize}
Our new \textit{subsampling online matrix factorization} algorithm (\somf) reduces the dimensionality of each of these steps, so that the
single-iteration complexity in $p$ depends on $q = \frac{p}{r}$ rather than
$p$. $r > 1$ is a \textit{reduction factor} that is close to the computational
speed-up per iteration in the large dimensional regime $p \gg k$.
Formally,
we randomly draw, at iteration $t$, a mask $\M_t$ that ``selects'' a random
subset of $\x_t$. $\M_t$ is a $\RR^{p\times p}$ random diagonal matrix, such
that each coefficient is a Bernouilli variable with parameter $\frac{1}{r}$,
normalized to be $1$ in expectation. With this definition at hand, $\M_t \x_t$
constitutes a non-biased, low-dimensional estimator of $\x_t$:
$
    \EE[\Vert \M_t \x_t \Vert_0] = \frac{p}{r} = q$, and $
    \EE[\M_t \x_t] = \x_t
$,
with $\Vert \cdot \Vert_0$ counting the number of non-zero coefficients. Thus, $r$ is the average proportion of observed features at each iteration.
We further define the pair of orthogonal projectors $\P_t \in \RR^{q \times p}$ and
$\P_t^\perp \in \RR^{(p - q)\times p}$ that projects $\RR^p$ onto $\mathrm{Im}(\M_t)$ and $\mathrm{Ker}(\M_t)$, which we will use for the dictionary update step.

In brief, \somf, defined in Alg.~\ref{alg:somf},
follows the outer loop of online matrix factorization, with the following major modifications at iteration $t$:
\begin{itemize}
    \item it uses $\M_t \x_t$ and low-size statistics instead of $\x_t$ to
    estimate the code $\balpha_t$ and the surrogate~$g_t$,
    \item it updates a subset of the dictionary $\P_t \D_{t-1}$ to reduce
    the surrogate value $\bar g_t(\D)$. Relevant parameters of $\bar g_t$ are computed
    using $\P_t \x_t$ and $\balpha_t$ only.
\end{itemize}

We describe in detail the new code computation and dictionary update steps. We then state convergence guarantees for \somf.
Those are non trivial to obtain as \somf is not an exact (stochastic) majorization-minimization
algorithm.

\begin{algorithm}[t]
    \begin{algorithmic}
        \Input Initial iterate $\D_0$, weight sequences ${(w_t)}_{t>0}$, ${(\gamma_c)}_{c>0}$,
        sample set ${\{\x^{(i)}\}}_{i> 0}$, \# iterations $T$.
        \For{$t$ from $1$ to $T$}
        \State Draw $\x_t = \x^{(i)}$ at random and $\M_t$ (see text)
        \State Update the regression parameters for sample $i$: $c^{(i)} \gets c^{(i)} + 1$, $\gamma \gets \gamma_{c^{(i)}}$
            \begin{equation*}
            (\bbeta_t^{(i)}, \G_t^{(i)}) \gets (1 - \gamma) (\bbeta_{t-1}^{(i)}, \G_{t-1}^{(i)})
             + \gamma (\D_{t-1}^\top \M_t \x^{(i)}, \D_{t-1}^\top \M_t \D_{t-1}),\,
             ( \bbeta_t, \G_t) \gets (\bar \bbeta_t^{(i)}, \bar \G_t^{(i)})
            \end{equation*}
        \State Compute the approximate code for $\x_t$: $\balpha_t \gets \argmin_{\balpha \in \RR^k}
            \frac{1}{2} \balpha^\top  \G_t \balpha -
            \balpha^\top  \bbeta_t + \lambda \, \Omega(\balpha).$
        \State Update the parameters of the aggregated surrogate $\bar g_t$:
        \begin{equation}
            \label{eq:somf_partial}
             \bar \C_t \gets (1 - w_t) \bar \C_t + w_t  \balpha_t  \balpha_t^\top. \qquad
             \P_t \bar \B_t \gets (1 - w_t ) \P_t \bar \B_t + w_t \P_t \x_t  \balpha_t^\top.
        \end{equation}
        \State Compute simultaneously (using (using \cite[Alg 2]{mensch_dictionary_2016}) for \nth{1} expression):
        \begin{equation}
            \label{eq:somf_minimization}
            \P_t \D_t \gets \argmin_{\D^r \in \mathcal{C}^r}
            \frac{1}{2} \trace ({\D^r}^\top (\D^r \bar \C_t - \bar \B_t)),\:
             \P_t^\perp \bar \B_t \gets (1 - w_t ) \P_t^\perp \bar \B_{t-1} + w_t \P_t^\perp \x_t  \balpha_t^\top.
        \end{equation}
        \EndFor
        \Output Final iterate $\D_T$.
    \end{algorithmic}
    \caption{Subsampled online matrix factorization (\somf)}\label{alg:somf}
\end{algorithm}
\paragraph{Code computation.}

In the online algorithm, $\balpha_t$ is obtained solving the linear regression problem
\begin{equation}
    \label{eq:regression}
    \balpha_t = \argmin_{\balpha \in \RR^k} \frac{1}{2} \balpha^\top \G_t^\star \balpha - \balpha^\top
\bbeta_t^\star + \lambda \Omega(\balpha),\quad\text{where}\quad
    \G_t^\star = \D_{t-1}^\top \D_{t-1}\:\text{and}\:\bbeta_t^\star = \D_{t-1}^\top \x_t
\end{equation}
For large $p$, computing $\G_t^\star$ and $\bbeta_t^\star$ dominates the complexity of the
code computation step. To reduce this complexity, we introduce \textit{estimators} for $\G_t$ and
$\bbeta_t$, computable at a cost proportional to $q$, whose use does not break convergence.
Recall that the sample $\x_t$ is drawn from a finite set of samples ${\{\x^{(i)}\}}_i$.
We estimate $\G_t^\star$ and $\bbeta_t^\star$ from $\M_t \x_t$ and data from previous iterations
 $s < t$ for which $\x^{(i)}$ was drawn. Namely, we keep in
memory $2n$ estimators, written ${( \G_t^{(i)}, \bbeta^{(i)}_t)}_{1\leq i \leq
n}$, observe the sample $i = i_t$ at iteration~$t$ and use it to
update the $i$-th estimators $\bar \G_t^{(i)}$, $\bar \bbeta^{(i)}_t$ following
\begin{equation}
        \bbeta_t^{(i)} = (1 - \gamma)  \G_{t-1}^{(i)} + \gamma \D_{t-1}^\top \M_t \x^{(i)},\qquad
        \G_t^{(i)} = (1 - \gamma)  \G_{t-1}^{(i)} + \gamma \D_{t-1}^\top \M_t \D_t^{(i)},
\end{equation}
where $\gamma$ is a weight
factor determined by the number of time sample $i$ has been previously observed at time
$t$. Precisely, given ${(\gamma_c)}_c$ a decreasing sequence of weights,
$\gamma \triangleq \gamma_{c^{(i)}_t}$, where
$
    c^{(i)}_t =
         | \lbrace s \leq t, \x_s = \x^{(i)} \rbrace |
$
All others estimators $\{\G^{(j)}_t, \bbeta^{(j)}_t\}_{j \neq i}$ are left unchanged from iteration
$t-1$. The set ${\{ \G_t^{(i)}, \bbeta^{(i)}_t\}}_{1\leq i\leq n}$ is used to define
the \textit{averaged} estimators at iteration $t$, related to sample $i$:
\begin{equation}
    \label{eq:agg-estimates}
     \G_t \triangleq  \G_t^{(i)} = \sum_{s \leq t, \x_s = \x^{(i)}} \gamma_{s,t}^{(i)} \D_{s-1}^\top \M_s \D_{s-1},\quad
     \bbeta_t \triangleq  \bbeta_t^{(i)} = \sum_{\substack{s \leq t, \x_s = \x^{(i)}}} \gamma_{s,t}^{(i)} \D_{s-1}^\top \M_s \x^{(i)},
\end{equation}
where
$\gamma_{s,t}^{(i)} = \gamma_{c^{(i)}_t} \prod_{s < t, \x_s = \x^{(i)}} (1 -
\gamma_{c^{(i)}_s})$. Replacing $(\G_t^\star, \bbeta_t^\star)$ by $(\G_t, \bbeta_t)$ in~\eqref{eq:regression}, $\balpha_t$
minimizes the masked loss averaged over the previous iterations where sample
$i$ appeared:
\begin{equation}
    \label{eq:approx-regression}
     \min_{\balpha \in \RR^k} \sum_{s\leq t\\
     \x_s = \x^{(i)}} \frac{\gamma_{s,t}^{(i)}}{2}
     \Vert \M_s(\x^{(i)} - \D_{s-1}^\top \balpha) \Vert_2^2
      + \lambda \Omega( \balpha ).
\end{equation}
The sequences ${(\G_t)}_t$ and ${(\bbeta_t)}_t$ are \textit{consistent}
estimations of ${(\G_t^\star)}_t$ and ${(\bbeta_t^\star)}_t$ --- consistency
arises from the fact that a single sample $\x^{(i)}$ is observed with different
masks along iterations. This was not the case
in the algorithm proposed in \cite{mensch_dictionary_2016}, which use estimators that does involve averaging from past data. Solving~\eqref{eq:approx-regression} is made closer and
closer to solving~\eqref{eq:regression}, in a manner that ensures the
correctness of the algorithm. Yet, computing the
estimators~\eqref{eq:agg-estimates} is $r$ times as costly as computing $\G_t^\star$ and $\bbeta_t^\star$
from~\eqref{eq:regression} and permits to speed up the code computation steop
close to $r$ times. The weight sequences $(w_t)_t$ and $(\gamma_c)_c$ are selected
 appropriately to ensure convergence.
For instance, we can set $w_t = \frac{1}{t^v}, \gamma_c = \frac{1}{c^{2.5 - 2v}}$,
with $v \in (\frac{3}{4}, 1)$.

\paragraph{Dictionary update.}In the original online algorithm, the whole dictionnary $\D_{t-1}$ is updated
at iteration $t$. To reduce the time complexity of this step,
we add a ``freezing'' constraint to the minimization of the quadratic function~\eqref{eq:full_quadratic}. Every
row $r$ of $\D$ that corresponds to an
\textit{unseen} row $r$ at iteration $r$ (such that $\M_t[r, r] = 0$) remains unchanged.
$\D_t$ is obtained thus solving
\begin{equation}
    \label{eq:dict-update-cons}
    \D_t \in \argmin_{\substack{\D \in \mathcal{C}\\
    \P_t^\perp \D = \P_t^\perp \D_{t-1}}}
     \frac{1}{2} \trace (\D^\top \D \bar \C_t)
      - \trace (\D^\top \bar \B_t),\,\text{with } \P_t\text{ orth. projector on }\mathrm{Im}(\M_t)
\end{equation}
With elastic-net ball constraints, solving \eqref{eq:dict-update-cons} reduces
to performing the partial dictionary update~\eqref{eq:somf_minimization} (Alg.~\ref{alg:somf}),
 with $\mathcal{C}^r = \{\D^r \in \RR^{k \times q}, \Vert {(\d^r)}^{(j)} \Vert \leq 1 -
\Vert \d^{(j)}_{t-1} \Vert + \Vert \P_t \d_{t-1}^{(j)} \Vert\}$.
We perform this update using a single pass of projected block coordinate descent with
blocks in the reduced space $\RR^{q}$. The dictionary update step is thus performed $r$ times faster than
the original algorithm.

\paragraph{Surrogate computation}The gradient we use to solve~\eqref{eq:somf_minimization} requires
to know only $\bar \C_t$ and $\P_t \bar \B_t$. We thus parallelize the partial update of the dictionary and
the update of $\P_t^\perp \bar \B_t$, using a second thread. The update of $\P_t \bar \B_t$
is performed in the main thread at a cost proportional to $q$. As the parallel computation
is dominated by dictionary update, this is enough to effectively reduce the computation time
of $\bar \B_t$ computation.

\paragraph{Convergence guarantees}

All in all, the three steps whose complexity depends on $p$ in the original algorithm now depends on $q$, which speeds-up a single iteration by a factor close to $r$. Yet \somf retain convergence guarantees.
We assume that there exists $\nu$ such that for all $t > 0$, $\D_t^\top \D_t \succ \nu \I$ (met in practice or by adding a small ridge regularization to~\eqref{eq:empirical-risk}), and make a technical data-independent hypothesis
on $(w_t)_t$ and $(\gamma_c)_c$ decay. Then \somf iterates converge in the same sense as in the original algorithm~\cite{mairal_online_2010}:
$\D_t \to \D_\infty \in \RR^{p\times k}$
and $\nabla \bar f(\D_\infty, \D - \D_\infty) \geq 0$
for any $\D \in \mathcal{C}$, where $\bar f$ is the empirical risk defined
in~\eqref{eq:empirical-risk}.
\begin{proof}[Proof sketch]
We use the aggregated nature of $\bar g_t$ and the fact that we can obtain a
\textit{geometric} rate of convergence for a single pass of projected block
coordinate descent \cite[\textit{e.g} see][]{richtarik_iteration_2014} to control the
terms $\D_t - \D_{t-1}$ and $\bar g_t(\D_t) - \bar g_t(\D_t^*)$, where $\D_t^* =
\argmin_{\D \in \mathcal{C}} \bar g_t$. We obtain $\bar g_t(\D_t) - \bar
g_t(\D_{t-1}) = \mathcal{O}(w_t)$, a crucial result on estimate stability. Simultaneously,
we show that the partial minimization yields the same result
as the full minimization for $t \to \infty$, as $\theta_t - \theta_t^\star \to 0$.
Using this result, with appropriate selection of $(w_t)_t$ and
$(\gamma_t)_t$, the noise induced by the use of estimators in
\eqref{eq:regression} can be bounded in the derivations of
\cite{mairal_online_2010}. We then write the difference
\eqref{eq:approx-regression} - \eqref{eq:regression} as the sum of a
\textit{lag} term and an empirical mean over ${\{\M_s\}}_{\x_s =\x_t}$. Both can be
bounded with appropriate selection of weights. Formally, ${(\bar g_t)}_t$ are
no longer upper-bounds of ${(\bar f_t)}_t$, but become so for $t \to \infty$, at a
sufficient rate to guarantee convergence.
\end{proof}

\section{Experiments}
\paragraph{Hyperspectral images.}We benchmark our algorithm by performing dictionary learning
on a large hyperspectral image. Dictionary learning is indeed used on
patches of hyperspectral images, as in
\cite{maggioni_nonlocal_2013}. Extracting $16{\times}16$ patches from an 1GB
hyperspectral image from the AVIRIS project with $224$ channels
yields samples of dimension $p = 57,000$. Figure 1 demonstrates that the newly
proposed algorithm is faster than the original, non-subsampled algorithm from
\cite{mairal_online_2010} by a factor close to $r = 4$. This speed-up is helped by the
redundancy in the different channels of the hyperspectral patches. In the first epochs, the proposed
method also outperforms the recent subsampled algorithm from
\cite{mensch_dictionary_2016}, thanks to the introduction of consistent
estimators in the code computation step. All algorithms are implemented in
Cython and benched on two cores. They cycle over $100,\!000$ normalized samples with minibatches of size~50. The
code for reproduction is available at \url{github.com/arthurmensch/modl}..

\begin{minipage}[c]{0.3\textwidth}
    \paragraph{Figure 1}
        Performing dictionary learning on hyperspectral data (224 channels, $16\times16$ patches) is faster with stochastic subsampling, and
        even faster with the newly proposed variance control.
    \label{fig:opt}
\end{minipage}%
\hfill%
\begin{minipage}[c]{0.5\textwidth}
    \includegraphics[width=\textwidth]{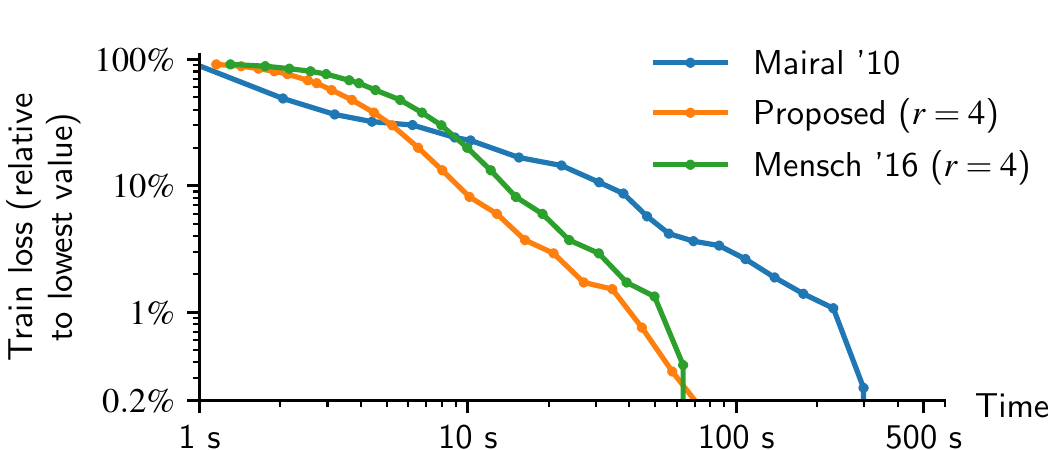}
\end{minipage}\hfill

\paragraph{Conclusion}The new \somf algorithm can efficiently factorize
large and tall matrices.
It preserves the speed gains of \cite{mensch_dictionary_2016} but has convergence guarantees that
are beneficial in practice.

\pagebreak

\bibliographystyle{abbrv}
\bibliography{opt_bib}

\begin{thebibliography}{10}

\bibitem{agarwal_learning_2014}
A.~Agarwal, A.~Anandkumar, P.~Jain, P.~Netrapalli, and R.~Tandon.
\newblock Learning sparsely used overcomplete dictionaries.
\newblock In {\em Conference on Learning Theory}, 2014.

\bibitem{defazio_saga:_2014}
A.~Defazio, F.~Bach, and S.~Lacoste-Julien.
\newblock Saga: A fast incremental gradient method with support for
  non-strongly convex composite objectives.
\newblock In {\em Advances in Neural Information Processing Systems}, 2014.

\bibitem{johnson_accelerating_2013}
R.~Johnson and T.~Zhang.
\newblock Accelerating stochastic gradient descent using predictive variance
  reduction.
\newblock In C.~J.~C. Burges, L.~Bottou, M.~Welling, Z.~Ghahramani, and K.~Q.
  Weinberger, editors, {\em Advances in {{Neural Information Processing
  Systems}}}. 2013.

\bibitem{maggioni_nonlocal_2013}
M.~Maggioni, V.~Katkovnik, K.~Egiazarian, and A.~Foi.
\newblock Nonlocal transform-domain filter for volumetric data denoising and
  reconstruction.
\newblock {\em IEEE Transactions on Image Processing}, 22(1):119--133, 2013.

\bibitem{mairal_stochastic_2013}
J.~Mairal.
\newblock Stochastic majorization-minimization algorithms for large-scale
  optimization.
\newblock In {\em Advances in {{Neural Information Processing Systems}}}, 2013.

\bibitem{mairal_online_2010}
J.~Mairal, F.~Bach, J.~Ponce, and G.~Sapiro.
\newblock Online learning for matrix factorization and sparse coding.
\newblock {\em The Journal of Machine Learning Research}, 11:19--60, 2010.

\bibitem{mensch_dictionary_2016}
A.~Mensch, J.~Mairal, B.~Thirion, and G.~Varoquaux.
\newblock Dictionary learning for massive matrix factorization.
\newblock In {\em {{International Conference}} on {{Machine Learning}}}, 2016.

\bibitem{olshausen_sparse_1997}
B.~A. Olshausen and D.~J. Field.
\newblock Sparse coding with an overcomplete basis set: {{A}} strategy employed
  by {{V1}}?
\newblock {\em Vision Research}, 37(23):3311--3325, 1997.

\bibitem{richtarik_iteration_2014}
P.~Richt{\'a}rik and M.~Tak{\'a}{\v c}.
\newblock Iteration complexity of randomized block-coordinate descent methods
  for minimizing a composite function.
\newblock {\em Mathematical Programming}, 144:1--38, 2014.

\bibitem{schmidt_minimizing_2013}
M.~Schmidt, N.~L. Roux, and F.~Bach.
\newblock Minimizing finite sums with the stochastic average gradient.
\newblock {\em arXiv preprint arXiv:1309.2388}, 2013.

\bibitem{tibshirani_regression_1996}
R.~Tibshirani.
\newblock Regression shrinkage and selection via the lasso.
\newblock {\em Journal of the Royal Statistical Society. Series B
  (Methodological)}, pages 267--288, 1996.

\bibitem{zou_sparse_2006}
H.~Zou, T.~Hastie, and R.~Tibshirani.
\newblock Sparse principal component analysis.
\newblock {\em Journal of computational and graphical statistics},
  15(2):265--286, 2006.

\end{thebibliography}
\end{document}